\title{Classification and Galois conjugacy of Hamming maps}

\author{Gareth A. Jones\\
School of Mathematics\\
University of Southampton\\
Southampton SO17  1BJ, U.K.\\
{\tt G.A.Jones@maths.soton.ac.uk}
}

\documentclass[12pt]{article}
\usepackage[latin1]{inputenc}
\usepackage{a4wide}
\usepackage{latexsym}
\usepackage{amsmath}
\usepackage{amssymb}
\usepackage{amsfonts}
\newtheorem{thm}{Theorem}[section]
\newtheorem{lemma}[thm]{Lemma}
\newtheorem{cor}[thm]{Corollary}
\newtheorem{prop}[thm]{Proposition}
\date{}
\begin{document} 

\maketitle

\begin{abstract}
\noindent We show that for each $d\geq 1$ the $d$-dimensional Hamming graph $H(d,q)$ has an orientably regular surface embedding if and only if $q$ is a prime power $p^e$. If $q>2$ there are  up to isomorphism $\phi(q-1)/e$ such maps, all constructed as Cayley maps for a $d$-dimensional vector space over the field $F_q$. We show that for each such pair $d, q$ the corresponding Bely\u{\i} pairs are conjugate under the action of the absolute Galois group ${\rm Gal}\,\overline{\bf Q}$, and we determine their minimal field of definition. We also classify the orientably regular embedding of merged Hamming graphs for $q>3$.
\end{abstract}


{\bf MSC classification:} Primary 20B25,  secondary 05C10, 05C25, 14H37, 14H55, 30F10.

{\bf Keywords:} Hamming graph, Hamming map, automorphism group, Galois group.

{\bf Running head:} Hamming maps\\

\section{Introduction}

A map $\cal M$ on an oriented surface is said to be {\sl orientably regular\/} if its orientation-preserving automorphism group ${\rm Aut}^+{\cal M}$ acts transitively on the arcs (directed edges) of $\cal M$. A standard problem in topological graph theory is that of classifying the orientably regular embeddings of a given class of arc-transitive graphs. This has been solved for several classes, for instance complete graphs $K_n$~\cite{Big, JJ}, complete bipartite graphs $K_{n,n}$~\cite{J2}, merged Johnson graphs $J(n,m)_I$~\cite{J1}, $n$-cubes $Q_n$ for $n$ odd~\cite{DKN}, and, according to a recent announcement~\cite{CCDKNW}, also $Q_n$ for $n$ even. Here we solve this problem for the $d$-dimensional Hamming graphs $H(d,q)$, where $q>2$ (when $q=2$ we have $H(d,2)\cong Q_d$ with a completely different classification): there are, up to isomorphism, $\phi(q-1)/e$ orientably regular embeddings when $q=p^e$ for some prime $p$, and none otherwise; as in the case $d=1$, when $H(1,q)\cong K_q$, these maps are all Cayley maps for a $d$-dimensional vector space over the field $F_q$ of $q$ elements. By contrast, Kwon~\cite{Kwo2} has recently classified the non-orientable regular embeddings of Hamming graphs, and in addition to the well-known non-orientable embeddings of $H(1,6)=K_6$ of type $\{3,5\}$ and $\{5,5\}$, obtained as antipodal quotients of the icosahedron and great dodecahedron~\cite{Jam}, there are two non-orientable embeddings of $H(2,6)$ of type $\{8,10\}$ and $\{10,10\}$.

We also consider the distance $k$ Hamming graphs $H(d,q)_k$ and the merged Hamming graphs $H(d,q)_K,\;\emptyset\neq K\subseteq D=\{1, 2, \ldots, d\}$: we show that if such a graph has an orientably regular embedding then $q$ is a prime power, and if $q\geq 4$ then the only orientably regular embeddings are those of $H(d,q)_1=H(d,q)$ described above, and the Biggs maps~\cite{Big} for the complete graph $H(d,q)_D=K_{q^d}$. The method of proof depends on a classification of the mergings of the Hamming association scheme by Muzychuk~\cite{Muz}, which applies only for $q\geq 4$; the case $q=3$ remains open.

According to Grothendieck's theory of {\it dessins d'enfants}~\cite{Gro, JS2}, maps on compact oriented surfaces correspond to algebraic curves defined over the field $\overline{\bf Q}$ of algebraic numbers; here we apply some general methods recently developed by Streit, Wolfart and the author~\cite{JSW}, based on Wilson's map operations~\cite{Wil}, to show that for $q>2$ the orientably regular embeddings of each $H(d,q)$ form an orbit under the absolute Galois group ${\rm Gal}\,\overline{\bf Q}$, and that their minimal field of definition is the splitting field of $p$ in the cyclotomic field of the $(q-1)$-th roots of $1$. In particular, they are defined over $\bf Q$ if and only if $q\leq 4$.

The author is grateful to J\"urgen Wolfart, and to the Mathematics Institute of the J.~W.~Goethe University, Frankfurt-am-Main, for hospitably arranging and funding a visit during which much of this research was carried out, and also to him and to Mikhail Klin and Manfred Streit for their helpful comments on early drafts of this paper. He also thanks the organisers of GEMS 09 for their invitation to announce these results at Tal\'e, Slovakia.

\section{Construction of Hamming maps}

The $d$-dimensional Hamming graph $H=H(d,q)$ has vertex set $V=Q^d$ where $Q$ is a set of $q$ elements ($q\geq 2$), with vertices $v=(v_i)$ and $w=(w_i)$ adjacent if and only if $v_i=w_i$ for all except one value of $i$. Each vertex thus has valency
$$n=d(q-1).$$
In any connected graph, the distance between two vertices $v$ and $w$ is defined to be the minimum number of edges in any path from $v$ to $w$; in this case it coincides with the {\sl Hamming distance}, the number of $i$ such that $v_i\neq w_i$. The automorphism group $A=A(d,q)$ of $H$ is the wreath product $S_q\wr S_d$ of the symmetric groups $S_q$ and $S_d$. This is a semidirect product $B\negthinspace :\negthinspace S_d$ of a normal subgroup $B=B_1\times\cdots\times B_d$, with the $i$-th direct factor $B_i\cong S_q$ acting naturally on $i$-th coordinates $v_i\in Q$ and fixing $j$-th coordinates $v_j$ for $j\neq i$, by a complement $S_d$ which permutes the coordinates $v_1,\ldots, v_d$ of each $v\in V$. If we take $Q$ to be a group, then $V$ is also a group, and $H$ is its Cayley graph with respect to the generating set $S$ consisting of the elements of $V$ with exactly one non-identity coordinate.

In certain cases one can form an orientably regular map which embeds $H$ in an oriented surface. Let $q$ be a prime power, and let $Q$ be the field $F=F_q$ of $q$ elements, so that the vertex set $V$ of $H$ is a $d$-dimensional vector space over $F$, and $S$ consists of the non-zero multiples of the standard basis vectors $e_1,\ldots, e_d$ of $V$. Let $\omega$ be a generator of the multiplicative group $F^*=F\setminus\{0\}\cong C_{q-1}$ of $F$, and let $M=M_{\omega}$ be the $d\times d$ monomial matrix over $F$ with non-zero entries $m_{i,i+1}=1$ for $i=1,\ldots, d-1$ and $m_{d,1}=\omega$. Then $M^d=\omega I$, so $M$ induces an invertible linear transformation of $V$ of order $n$, permuting the elements of $S$ in a single cycle. This cyclic ordering
\begin{equation}
e_1,\; e_1M=e_2,\;\ldots,\; e_1M^{d-1}=e_d,\; e_1M^d=\omega e_1,\; \ldots,\; e_1M^{d(q-1)-1}=\omega^{q-2}e_d
\end{equation}
of $S$ determines a Cayley map ${\cal H}={\cal H}(d,\omega)$ for $H$: the cyclic order of the neighbours of any vertex $v$, induced by the local orientation around $v$, is obtained by adding each vector in $(1)$ to $v$.  This map, which we will call a {\sl Hamming map}, is orientably regular since $M$ induces an automorphism of the group $V$ (see~\cite[Theorem 16-27]{Whi}). The orientation-preserving automorphism group $G={\rm Aut}^+{\cal H}(d,\omega)$ of ${\cal H}$ is a semidirect product $V\negthinspace :\negthinspace G_{\bf 0}$ of an elementary abelian normal subgroup $V\cong F^d\cong (C_p)^{de}$, acting regularly by translations on the vertices, by a complement $G_{\bf 0}\cong C_n$ fixing the vertex ${\bf 0}=(0,\ldots,0)$; a generator of $G_{\bf 0}$ permutes the neighbours of $\bf 0$ in a single cycle, given by $(1)$.

These Hamming maps are generalisations of the orientably regular embeddings ${\cal H}(1,\omega)$ of the complete graphs $K_q\cong H(1,q)$ introduced by Biggs in~\cite{Big}. James and the author showed in~\cite{JJ} that the maps ${\cal H}(1,\omega)$ are the only orientably regular embeddings of complete graphs, and here (for $q>2$) we extend this result to all dimensions $d\geq 1$:

\begin{thm}
Let $\cal M$ be an orientably regular embedding of the $d$-dimensional Hamming graph $H(d,q)$. Then $q$ is a prime power, and if $q>2$ then $\cal M$ is isomorphic to a Hamming map ${\cal H}(d,\omega)$ where $\omega$ is a generator of the multiplicative group $F_q^*$.
\end{thm}

The proof, given in \S 6, uses finite group theory, including basic properties of Frobenius groups, though it is independent of the classification of finite simple groups. This result has been obtained independently by Kwon~\cite{Kwo} for odd $q$, using a different method which represents oriented maps as pairs of permutations. The condition $q>2$ is necessary here, at least for $d\geq 3$, since in this case there are more orientably regular embeddings of the $d$-dimensional hypercube $Q_d\cong H(d,2)$ than the single Hamming map ${\cal H}(d,\omega)$, with $\omega=1\in F_2$: for instance if $d$ is odd then as shown by Du, Kwak and Nedela~\cite{DKN} there are, up to isomorphism, $2^r$ orientably regular embeddings of $Q_d$, where $r$ is the number of distinct primes dividing $d$. The freedom to create additional embeddings when $q=2$ seems to depend on $H(d,2)$ being bipartite, which is not the case for $H(d,q)$ when $q>2$.

For any prime power $q$ the number of distinct Hamming maps ${\cal H}(d,\omega)$ is $\phi(q-1)$, the number of generators $\omega$ of $F_q^*$, where $\phi$ is Euler's function. If $\omega'=\omega^{\gamma}$ for some field automorphism $\gamma$ of $F_q$, then letting $\gamma$ act naturally on $V$ induces an isomorphism ${\cal H}(d,\omega)\to{\cal H}(d,\omega')$. Our second main theorem, proved in \S 7 and also generalising a result for $d=1$ in~\cite{JJ}, shows that these are the only cases in which ${\cal H}(d,\omega)\cong{\cal H}(d,\omega')$:

\begin{thm} Let $\omega$ and $\omega'$ be generators of the multiplicative group of the field $F_q$. Then ${\cal H}(d,\omega)\cong{\cal H}(d,\omega')$ if and only if $\omega$ and $\omega'$ are conjugate under the Galois group of $F_q$.
\end{thm}

 If $q=p^e$ for some prime $p$ then the Galois group ${\rm Gal}\,F_q$ is cyclic of order $e$, generated by the Frobenius automorphism $t\mapsto t^p$. Since this group acts fixed-point-freely on the generators of $F_q^*$, we have the following immediate corollary to Theorems~2.1 and 2.2:

\begin{cor}
If $q=p^e>2$ for some prime $p$ then there are, up to isomorphism, exactly $\phi(q-1)/e$ orientably regular embeddings of $H(d,q)$ for each $d\geq 1$.
\end{cor}

It is, perhaps, a little surprising that this number is independent of the dimension $d$. Like Corollary~2.3, the following result, proved in \S 3, generalises one proved for $d=1$ in~\cite{JJ}. An orientably regular map $\cal M$ is {\sl reflexible\/} if it has an orientation-reversing automorphism. 

\begin{cor}
For each $d\geq 1$, the only reflexible Hamming maps are the unique Hamming maps for $q=2, 3$ or $4$. In particular, the only reflexible embeddings of $H(d,q)$ for $q>2$ are the Hamming maps with $q=3$ or $4$.
\end{cor}

\section{Properties of Hamming maps}

In any orientably regular map $\cal M$, the faces have the same valency $m$, the vertices have the same valency $n$, and the Petrie polygons (closed zig-zag paths) have the same length $l$. We then say that $\cal M$ has {\sl type\/} $\{m,n\}_l$, or  simply $\{m,n\}$. We now compute the type of the map ${\cal H}(d,\omega)$. 

\begin{lemma}
Let $\cal M$ be a Cayley map for an abelian group $A$, in which the cyclic order of the generators is given by successive iterates $\alpha^i(s)$ of an automorphism $\alpha$ of $A$ applied to an element $s\in A$. Then
\begin{description}
\item[(a)] if the automorphism $-\alpha: a\mapsto -\alpha(a)$ of $A$ fixes only $0$ in $A$ then the face valency $m$ is the order of this automorphism;
\item[(b)] the Petrie length $l$ is twice the order of the element $(\alpha-1)(s)=\alpha(s)-s$ in $A$.
\end{description}
\end{lemma}

\noindent{\sl Proof.} By its construction, $\cal M$ is orientably regular (see~\cite[Theorem 16-27]{Whi}), so all faces have the same valency $m$. Successive vertices around one particular face are
\[s,\; 0,\; \alpha(s),\; \alpha(s)-\alpha^2(s),\; \alpha(s)-\alpha^2(s)+\alpha^3(s),\; \ldots\]
so $m$ is the least $j>0$ such that
\begin{equation}
\alpha(s)-\alpha^2(s)+\alpha^3(s)-\cdots-(-1)^j\alpha^j(s)=0.
\end{equation}
If $(2)$ holds, then applying $\alpha^{-1}$ and adding, we see that $(-\alpha)^j(s)=s$; since the automorphism $-\alpha$ commutes with $\alpha$ this implies that $(-\alpha)^j(\alpha^i(s))=\alpha^i(s)$ for all $i$, so $(-\alpha)^j=1$ since the elements $\alpha^i(s)$ generate $A$. Conversely, if $(-\alpha)^j=1$ then
\[(\alpha^{-1}+1)(\alpha(s)-\alpha^2(s)+\alpha^3(s)-\cdots-(-1)^j\alpha^j(s))=s-(-1)^j\alpha^j(s)=0,\]
so $\alpha^{-1}$, and hence also $\alpha$, inverts the left-hand side of $(2)$. It follows that if $\alpha$ inverts only $0\in A$ then $(2)$ holds, so $m$ is the order of $-\alpha$.

Similarly all Petrie polygons in $\cal M$ have the same length $l$. A typical Petrie polygon has successive vertices
\[s,\; 0,\; \alpha(s),\; \alpha(s)-s,\; 2\alpha(s)-s,\; 2\alpha(s)-2s,\; \ldots,\; j\alpha(s)-(j-1)s,\; j\alpha(s)-js,\;\ldots,\]
so $l$ is twice the order of $(\alpha-1)(s)$. \hfill$\square$

\medskip

The vertices of the Hamming map ${\cal H}(d,\omega)$ have valency $n=d(q-1)$, and Lemma~3.1(b) shows that the Petrie length is $l=2p$, where $q=p^e$ and $p$ is prime.

Lemma~3.1(a) applies to the Hamming maps ${\cal H}(d,\omega)$ for $q>3$, and also for $q=3$ if $d$ is even: since the matrix $M=M_{\omega}$ has characteristic polynomial $\lambda^d-\omega$ we see that $-1$ cannot be an eigenvalue in these cases. Now $(-M)^d=(-1)^d\omega I$, and this is $\omega I$ if $q$ or $d$ is even, giving $m=n$; if $q$ and $d$ are both odd then $(-M)^d=-\omega I$ with $\omega^{(q-1)/2}=-1$, so $m=n$ or $n/2$ as $q\equiv 1$ or $-1$ mod~$(4)$. Thus $m=n$ in all these cases except when $d$ is odd and $q\equiv -1$ mod~$(4)$, in which case $m=n/2$.

To deal with the exceptional cases, first let $q=3$ with $d$ odd. Then $\omega=-1$ so $M^d=-I$ and hence $(-M)^d=I$; directly calculating the left-hand side of $(2)$ shows that $m=3d$. Finally, if $q=2$ then $\omega=1$, and this time calculating the left-hand side of $(2)$ gives $m=2d$. To summarise, $m=n$ unless
\begin{itemize}
\item $d$ is odd and $3<q\equiv -1$ mod~$(4)$, in which case $m=n/2$, or
\item $d$ is odd and $q=3$, in which case $m=3d$, or
\item $q=2$, in which case $m=2d$.
\end{itemize}

Knowing its type, we can now compute the Euler characteristic and hence the genus of ${\cal H}(d,\omega)$. The numbers of vertices, edges and faces are $q^d$, $d(q-1)q^d/2$ and $d(q-1)q^d/m$ respectively. It follows that if $m=d(q-1)$ then ${\cal H}(d,\omega)$ has characteristic
\[\chi=2q^d-\frac{d(q-1)q^d}{2}=\frac{q^d}{2}\bigl(4-d(q-1)\bigr)\]
and hence genus
\[g=1+\frac{q^d}{4}\bigl(d(q-1)-4\bigr).\]
If  $3<q\equiv -1$ mod~$(4)$ and $d$ is odd then $m=d(q-1)/2$, so ${\cal H}(d,\omega)$ has characteristic
\[\chi=3q^d-\frac{d(q-1)q^d}{2}=\frac{q^d}{2}\bigl(6-d(q-1)\bigr)\]
and hence genus
\[g=1+\frac{q^d}{4}\bigl(d(q-1)-6\bigr).\]
If  $q=3$ and $d$ is odd then $m=3d$, so
\[\chi=3^d-3^dd+2.3^{d-1}=3^{d-1}(5-3d)\]
giving
\[g=1+\frac{3^{d-1}(3d-5)}{2}.\]
If $q=2$ then $m=2d$, so
\[\chi=2^d-2^{d-1}d+2^{d-1}=2^{d-1}(3-d)\]
and hence
\[g=1+2^{d-2}(d-3).\]

\medskip

We now consider the effect of Wilson's operations $H_j$, which act on maps by preserving the graph and raising the cyclic order of neighbours of each vertex to its $j$-th power, where $j$ is coprime to the valency~\cite{Wil}.

\begin{lemma}
$H_j({\cal H}(d,\omega))\cong{\cal H}(d,\omega^j)$ for all $j$ coprime to $d(q-1)$.
\end{lemma}

\noindent{\sl Proof.} By raising the cyclic order of neighbours of $0$ to its $j$-th power we do the same to the induced $d$-cycle on the coordinate places, resulting in another $d$-cycle since $j$ is coprime to $d$, and we replace $\omega$ with $\omega^j$, which is a generator of $F_q^*$ since $j$ is coprime to $q-1$. This gives the required isomorphism. \hfill$\square$

\begin{cor} The mirror image of ${\cal H}(d,\omega)$ is isomorphic to ${\cal H}(d,\omega^{-1})$.  \hfill$\square$
\end{cor}

If $q=2$ or $3$ then $\omega=\pm 1$, so this implies that ${\cal H}(d,\omega)$ is reflexible; this also applies when $q=4$ since $\omega^{-1}=\omega^2$ is conjugate to $\omega$ under the Galois group of $F_4$ (see \S 2). If $q\geq 5$ then $\omega^{-1}$ and $\omega$ are never conjugate, so in this case Theorem~2.2 implies that ${\cal H}(d,\omega)$ and its mirror image are non-isomorphic, forming a chiral pair. This, together with Theorem~2.1, proves Corollary~2.4.

\begin{lemma}
For a given pair $d$ and $q$, the automorphism groups of the Hamming maps ${\cal H}(d,\omega)$ are mutually isomorphic.
\end{lemma}

\noindent{\sl Proof.} Let ${\cal H}(d,\omega)$ and ${\cal H}(d,\omega')$ be Hamming maps for a given pair $d$ and $q$. Their automorphism groups $G$ and $G'$ are semidirect products of $V=F_q^d$ by $C_n$, with a generator of $C_n$ acting linearly on $V$ as the matrix $M_{\omega}$ or $M_{\omega'}$ with respect to the standard basis $e_1,\ldots, e_d$.
We have $\omega'=\omega^j$ for some $j$ coprime to $q-1$, and we can choose $j$ also to be coprime to $d$. Then $M_{\omega}^j$ acts linearly on $V$ as the matrix $M_{\omega'}$ with respect to the basis $e_1,\, e_1M_{\omega}^j,\, e_1M_{\omega}^{2j},\ldots,\, e_1M_{\omega}^{(d-1)j}$, so applying the appropriate change-of-basis automorphism of $V$ and raising the elements of $G_0$ to their $j$-th powers gives the required isomorphism $G\to G'$.\hfill$\square$

\medskip

We will denote this common automorphism group by $G=G(d,q)$.

\section{Examples}

\noindent{\bf (a)} First we will briefly describe the maps ${\cal H}(d,\omega)$ of genus $g$ up to $101$, to allow comparison with Conder's classification~\cite{Con} of orientably regular maps in the range $2\leq g\leq 101$; the cases $g=0$ and $1$ are well-known, and can be found in~\cite{CM}. We exclude the Hamming maps with $d=1$ since these are the well-known complete maps discovered by Biggs~\cite{Big, JJ}. In most cases, a few numerical parameters such as the genus, type (including Petrie length), and number of automorphisms are sufficient to identify the map in Conder's lists, but occasionally one also needs to use the presentations of the automorphism groups in~\cite{Con}. These examples are presented in increasing order of $d$.

\medskip

\noindent If $q=2$ the unique Hamming map ${\cal H}(2,\omega)$ has type $\{4,2\}_4$ and genus $0$ with $|G|=8$; this is the reflexible map $\{4,2\}$ in~\cite[Ch.~8]{CM}, the spherical embedding of a cycle of length $4$.

\smallskip

\noindent If $q=3$ the unique Hamming map ${\cal H}(2,\omega)$ has type $\{4,4\}_6$ and genus $1$ with $|G|=36$; this is the reflexible map $\{4,4\}_{3,0}$ in~\cite[Ch.~8]{CM}, also a Paley map~\cite[\S 16.8]{Whi} for the Paley graph $P_9$. 

\smallskip

\noindent  If $q=4$ the unique Hamming map ${\cal H}(2,\omega)$ has type $\{6,6\}_4$ and genus $9$ with $|G|=96$; this is the self-dual map R9.18 in Conder's list of reflexible maps~\cite{Con}.

\smallskip

\noindent  If $q=5$ the two Hamming maps ${\cal H}(2,\omega)$ have type $\{8,8\}_{10}$ and genus $26$ with $|G|=200$; these form the chiral pair C26.1 in Conder's list of chiral maps~\cite{Con}.

\smallskip

\noindent  If $q=7$ the two Hamming maps ${\cal H}(2,\omega)$ have type $\{12,12\}_{14}$ and genus $99$ with $|G|=588$; these form the chiral pair C99.2 in~\cite{Con}.

\smallskip

\noindent  If $q=2$ the unique Hamming map ${\cal H}(3,\omega)$ has type $\{6,3\}_4$ and genus $1$ with $|G|=24$; this is the reflexible map $\{6,3\}_{2,0}$ in~\cite[Ch.~8]{CM}, the Petrie dual of the cube $\{4,3\}$.

\smallskip

\noindent  If $q=3$ the unique Hamming map ${\cal H}(3,\omega)$ has type $\{9,6\}_6$ and genus $19$ with $|G|=162$; this is the dual of the reflexible map R19.15 in~\cite{Con}.

\smallskip

\noindent If $q=4$ the unique Hamming map ${\cal H}(3,\omega)$ has type $\{9,9\}_4$ and genus $81$ with $|G|=576$; this is the self-dual reflexible map R81.125 in~\cite{Con}.

\smallskip

\noindent If $q=2$ the unique Hamming map ${\cal H}(4,\omega)$ has type $\{8,4\}_4$ and genus $5$ with $|G|=64$; this map, an embedding of the hypercube $Q_4$, is the dual of the reflexible map R5.5 in~\cite{Con}.

\smallskip

\noindent  If $q=3$ the unique Hamming map ${\cal H}(4,\omega)$ has type $\{8,8\}_6$ and genus $82$ with $|G|=648$; this is the self-dual reflexible map R82.50 in~\cite{Con}.

\smallskip

\noindent If $q=2$ the unique Hamming map ${\cal H}(5,\omega)$ has type $\{10,5\}_4$ and genus $17$ with $|G|=160$; this map, an embedding of $Q_5$, is the dual of the reflexible map R17.18 in~\cite{Con}.

\smallskip

\noindent  If $q=2$ the unique Hamming map ${\cal H}(6,\omega)$ has type $\{12,6\}_4$ and genus $49$ with $|G|=384$; this map, an embedding of $Q_6$, is the dual of the reflexible map R49.39 in~\cite{Con}.

\medskip

\noindent{\bf(b)} To see a class of examples for which there is more than one chiral pair of embeddings, let $q=25$. For each $d\geq 1$ there are $\phi(24)/2=4$ non-isomorphic Hamming maps ${\cal H}(d,\omega)$. For instance, the embeddings ${\cal H}(1,\omega)$ of the complete graph $K_{25}$ considered in~\cite[\S 8]{JJ} have type $\{24,24\}_{10}$ and genus $126$, while the maps ${\cal H}(2,\omega)$ have type $\{48,48\}_{10}$ and genus $6876$. For each $d$ these four maps correspond to the four orbits of ${\rm Gal}\,F_{25}\cong C_2$ on the $\phi(24)=8$ generators of $F_{25}^*$, or equivalently to the four irreducible quadratic factors $t^2\pm t+2$ and $t^2\pm 2t-2$ of the cyclotomic polynomial $\Phi_{24}(t)=t^8-t^4+1$ over $F_5$, each arising as the minimal polynomial of a conjugate pair of generators. The maps form two chiral pairs, each pair corresponding to a mutually inverse pair of generators of $F_{25}^*$, and transposed with the other pair by Wilson's operation $H_7$.

\section{Preliminaries for the classification}

For each $i=1,\ldots, d$ let $\sim_i$ be the equivalence relation on the vertex set $V$ of $H=H(d,q)$ defined by $v\sim_i w$ if and only if $v_i=w_i$, and let $\pi_i$ be the corresponding partition of $V$, consisting of $q$ equivalence classes $[v]_i$ of size $q^{d-1}$. The automorphism group $A={\rm Aut}\,H$ permutes these $d$ partitions $\pi_i$, acting as the symmetric group $S_d$; the kernel of its action is the normal subgroup $B=B_1\times\cdots\times B_d$ of $A$, with $B_i$ acting faithfully as $S_q$ on the $q$ classes in $\pi_i$, while the other factors $B_j$ leave each of these classes invariant.

For each $i=1,\ldots, d$ let $\approx_i$ be the conjunction (or intersection) of the equivalence relations $\sim_j$ for $j\neq i$; this is the equivalence relation on $V$ defined by $v\approx_i w$ if and only if $v_j=w_j$ for all $j\neq i$. The corresponding partition $\Pi_i$ of $V$ consists of $q^{d-1}$ equivalence classes $[[v]]_i$ of size $q$. Again, $A$ permutes these $d$ partitions $\Pi_i$ as $S_d$, with kernel $B$, but in this case $B_i$ is the kernel of the action of $B$ on the classes in $\Pi_i$. 

Each class $[[v]]_i$ in $\Pi_i$ is a maximal clique (complete subgraph) in $H$. For each $k=0,\ldots, d$ let $H_k(v)$ denote the set of vertices at distance $k$ from $v$. Then the set $H_1(v)$ of neighbours of $v$ in $H$ is the disjoint union of the cliques $[[v]]_i^*=[[v]]_i\setminus\{v\}$ for $i=1,\ldots, d$. (This fact easily implies the well-known result that ${\rm Aut}\,H\cong S_q\wr S_d$.)

In order to prove Theorem~2.1 we need the following result:

\begin{lemma}
The automorphisms of a Hamming graph $H$ fixing a vertex $v$ are represented faithfully on  $H_1(v)$.
\end{lemma}

\noindent{\sl Proof.} Suppose that a graph automorphism $g$ fixes $v$ and all its neighbours. For $k\geq 2$, each vertex in $H_k(v)$ is uniquely determined by its $k$ neighbours in $H_{k-1}(v)$, so it follows by induction on $k$ that $g$ fixes every vertex of $H$. \hfill$\square$

\section{Proof of Theorem 2.1}

Let $G={\rm Aut}^+{\cal M}$ for some orientably regular embedding $\cal M$ of $H(d,q)$, where $q>2$. As explained earlier, it is shown in~\cite{JJ} that Theorem~2.1 is true for $d=1$, so we may assume that $d\geq 2$. 

The subgroup $G_v$ of $G$ stabilising a vertex $v\in V$ is a cyclic group of order $n=d(q-1)$, acting regularly on the set $H_1(v)$ of neighbours of $v$. Since $G\leq A$, $G_v$ acts imprimitively on $H_1(v)$, permuting the $d$ cliques $[[v]]_i^*$ transitively. Let $K=G\cap B$, the kernel of the action of $G$ on the $d$ partitions $\Pi_i$ of $V$. Then $K_v=G_v\cap K$ is a cyclic group of order $q-1$ acting regularly on each of the $d$ sets $[[v]]_i^*$. It follows from Lemma~5.1 that $K_{vw}=1$ for each $w\in V\setminus\{v\}$. This will show that $K$ acts on $V$ as a Frobenius group, provided we can show that $K$ acts transitively but not regularly on $V$. If $w\approx_i v$ then $K_w$ acts regularly on $[[w]]_i^*=[[v]]_i\setminus\{w\}$; the subsets $[[v]]_i^*$ and $[[w]]_i^*$ of $[[v]]_i=[[w]]_i$ have non-empty intersection since $q>2$, so all elements of $[[v]]_i$ are in the same orbit of $K$; since the transitive closure of the $d$ equivalence relations $\approx_i$ is the universal relation, $K$ is transitive on $V$, and since $|K_v|=q-1>1$ it follows that $K$ acts as a Frobenius group on $V$.

As a Frobenius group, $K$ has a normal subgroup $N$ acting regularly on $V$, namely the Frobenius kernel~\cite[V.7.6, V.8.2]{Hup}. Since $N$ is a Hall subgroup of $K$, it is a characteristic subgroup of $K$~\cite[V.8.3]{Hup} and hence it is normal in $G$. It follows that $G$ is a semidirect product of $N$ by $G_v$ for each $v\in V$.

Let us choose a particular element of $Q$, denoted by $0_Q$, and let $0$ denote the vertex $(0_Q,\ldots,0_Q)\in V$. We can identify $V$ with $N$ so that $0$ is the identity element, $N$ acts on $V$ by right multiplication, and $G_0$ acts on $V$ by conjugation. For each $i$ the equivalence class $[[0]]_i$ containing $0$ is identified with a subgroup $N_i$ of order $q$, namely the subgroup of $N$ preserving $[[0]]_i$. Since $K_0$ acts transitively by conjugation on the non-identity elements of $N_i$, this subgroup is elementary abelian, so $q=p^e$ for some prime $p$. This proves the first part of Theorem~2.1.

In order to prove the second part, we first need to show that $N$ induces a regular permutation group on the set of $q$ classes in each partition $\pi_i$. There is a simple proof of this when $q$ is odd, using the fact that the order $q-1$ of the Frobenius complement $K_0$ is then even, so that the Frobenius kernel $N$ is abelian~\cite[V.8.18(a)]{Hup}. However, this argument fails when $q=2^e$, so instead we give an alternative argument which applies in all cases.

For each $i$, $K$ permutes the $q$ classes in the partition $\pi_i$ as a transitive permutation group $K^{(i)}$. The subgroup $K_0$ preserves the class $[0]_i$, and permutes the other $q-1$ classes regularly, so $K^{(i)}$ is a $2$-transitive group. The normal subgroup $N$ of $K$ has order $q^d=p^{de}$, so it induces a normal $p$-group $N^{(i)}\leq K^{(i)}$, and this contains the regular subgroup $N_i^{(i)}$ induced by the subgroup $N_i$ on $\pi_i$. If $N^{(i)}>N_i^{(i)}$ then the stabiliser of $[0]_i$ in $K^{(i)}$ contains a non-trivial normal $p$-subgroup; the orbits of this subgroup on the remaining $q-1$ classes must have the same length $l>1$ dividing $q-1$, which is impossible since $q-1$ is coprime to $p$. Thus $N^{(i)}=N_i^{(i)}$, so $N$ induces a regular permutation group $N^{(i)}$ on $\pi_i$. Each subgroup $N_j$ ($j\neq i$) of $N$ sends $0$ to elements of $[[0]]_j\subseteq[0]_i$, so it preserves the class $[0]_i$ and is therefore contained in the kernel of this action of $N$ on $\pi_i$. This is true for all $i\neq j$, so $N_j$ is contained in the kernel $N\cap B_j$ of the action of $N$ on $\Pi_j$. Thus $N_j\leq B_j$ for each $j$, so the subgroups $N_1,\ldots, N_d$ generate their direct product $N_1\times\cdots\times N_d$ in $N$, and comparing orders we see that $N=N_1\times\cdots\times N_d$. Thus $N$ is an elementary abelian $p$-group of rank $de$.

Let $x$ be the standard generator of $G_{\bf 0}$, permuting the edges of $\cal M$ incident with $0$ by following the local orientation around this vertex. Then $x$ permutes the sets $[[0]]_i^*$ in a cycle of length $d$, and by renumbering if necessary (equivalently, by applying an automorphism of $H$) we may assume that it acts on the subscripts $i$ as the cycle $(1, 2, \ldots, d)$.

The subgroup $K_0=\langle x^d\rangle$ of $K$, being cyclic and acting transitively on the non-identity elements of each $N_i\cong (C_p)^d$, must act on each $N_i$ as a Singer cycle~\cite[II.3.10, II.7.3]{Hup}, so in particular one can identify $N_1$ with the field $F=F_q$ so that $x^d$ acts on $N_1$ as multiplication by a generator $\omega$ of $F^*$. By composing this identification $F\to N_1$ with the bijection $N_1\to N_i$ induced by $x^{i-1}$ for each $i=2,\ldots, d$, we get an identification $F\to N_i$ of $N_i$ with $F$ so that $x^d$ also acts on $N_i$ as multiplication by $\omega$. Then $x:N_i\to N_{i+1}$ acts as the identity on $F$ for $i=1,\ldots, d-1$, while $x:N_d\to N_1$ acts as multiplication by $\omega$. By Lemma~5.1 this action of $x$ on $N_1\cup\cdots\cup N_d$ has a unique extension to a graph automorphism of $H$, namely the cyclic permutation $(1)$ defined in \S 2. This shows that the vertex-set $V=N$ can be regarded as the $d$-dimensional vector space $F^d$ over $F$, with $x$ acting as the matrix $M_{\omega}$ with respect to the standard basis, so ${\cal M}\cong{\cal H}(d,\omega)$. \hfill$\square$

\section{Proof of Theorem 2.2}

We need to show that if ${\cal H}(d,\omega)\cong{\cal H}(d,\omega')$ where $\omega$ and $\omega'$ are generators of $F^*$, then $\omega$ and $\omega'$ are in the same orbit of ${\rm Gal}\,F$. Equivalently, we need to show that a Hamming map $\cal M$ determines $\omega$ uniquely, up to automorphisms of $F$.

Given such a map $\cal M$, let us choose a pair of adjacent vertices, and label them with the elements $0$ and $e_1$ of $V$. (The particular choice is immaterial, since $\cal M$ is orientably regular.) In the proof of Theorem~2.1, the unique maximal clique $[[0]]_1$ of $H$ containing $0$ and $e_1$ is identified with the $1$-dimensional subspace $N_1$ of $N=V=F^d$ spanned by $e_1$; we can therefore identify it with the field $F=F_q$, by identifying each element $\lambda e_1\;(\lambda\in F)$ of $[[0]]_1$ with $\lambda$. We need to find this identification explicitly, in order to determine $\omega$.

By following the orientation of $\cal M$ around $0$ we obtain the neighbours of $0$ in the cyclic order given by $(1)$; the $d$-th power of this cyclic permutation, starting with the chosen vertex $e_1$, therefore identifies the elements of $[[0]]_1^*$ with the successive powers $1, \omega, \omega^2,\ldots, \omega^{q-2}$ of $\omega$. This gives us a cyclic group structure on $[[0]]_1^*$, which  we can regard as the multiplicative group $F^*$ of $F$.

In order to determine the additive structure of $[[0]]_1$ it is sufficient to find the difference of each distinct pair of elements $v=\omega^i$ and $w=\omega^j$ of $[[0]]_1^*$, expressing $v-w$ as a power of $\omega$. In the cyclic order of the neighbours of $w$ in $\cal M$, the vertex $0=w+(-w)$ must be followed, $dk$ terms later for some integer $k$, by the vertex $v$, so that $v=w+\omega^k(-w)$. Thus $v-w=\omega^k(-w)$, and this is either $\omega^{j+k}$ or $\omega^{j+k+(q-1)/2}$ as $q$ is even or odd, since $-1=\omega^{(q-1)/2}$ in the latter case. This shows that $\cal M$ uniquely determines the field structure of $F$, so it determines the minimal polynomial of $\omega$ over the prime field $F_p$, and hence it determines $\omega$ up to an automorphism of $F$. \hfill$\square$

\section{Characterisation by valency and automorphisms}

For future applications, we need to show that the Hamming maps ${\cal H}(d,\omega)$ are characterised among all orientably regular maps  by their valency and their orientation-preserving automorphism group. We first summarise some basic general facts about orientably regular maps; for background, see~\cite{JS1}, for instance.

For any group $G$, the orientably regular maps $\cal M$ with ${\rm Aut}^+{\cal M}\cong G$ correspond to the generating pairs $x, y$ for $G$ such that $y$ has order $2$. Here $x$ is a rotation fixing a vertex $v$ of $\cal M$, sending each incident edge to the next incident edge according to the local orientation around $v$, while $y$ is a half-turn, reversing one of these incident edges, so that $z=(xy)^{-1}$ is a rotation preserving an incident face. We will call $x$ and $y$ {\sl standard generators\/} of $G$. Conversely, given generators $x$ and $y$ of a group $G$ with $y^2=1$ one can construct a map $\cal M$, with arcs corresponding to the elements of $G$, and vertices, edges and faces corresponding to the cosets in $G$ of the cyclic subgroups generated by $x, y$ and $z$; this map has type $\{m,n\}$ where $m$ and $n$ are the orders of $z$ and $x$. Two such maps are isomorphic if and only if the corresponding pairs of standard generators are equivalent under ${\rm Aut}\,G$.

We will apply this theory to the common automorphism group $G=G(d,q)$ of the Hamming maps for $H(d,q)$, described in \S 3.

\begin{lemma} If $q>2$ then the matrix $M_{\omega}$ satisfies
\[\sum_{i=0}^{n-1}M_{\omega}^i=0.\]
\end{lemma}

\noindent{\sl Proof.} By direct calculation, each entry of the matrix $\sum_{i=0}^{n-1}M_{\omega}^i$ is equal to $\mu=\sum_{i=0}^{q-2}\omega^i$. Since $\omega^{q-1}=1$ we have $\omega\mu=\mu$, and hence $\mu=0$ since $\omega\neq 1$ when $q>2$. \hfill$\square$

\medskip

[If $q=2$ then $\sum_{i=0}^{n-1}M_{\omega}^i$ is the $d\times d$ matrix $J$ with all entries equal to $1$.]

\medskip

Each element of the semidirect product $G=V\negthinspace :\negthinspace G_0$ has a unique factorisation $vg$ where $v\in V$ and $g\in G_0$. We will use this to determine the elements of order $n$ and $2$ in $G$.

\begin{cor} An element $vg\in G$, with $v\in V$ and $g\in G_{\bf 0}$, has order $n$ if and only if $g$ has order $n$.
\end{cor}

\noindent{\sl Proof.} For each $k\geq 1$ we have
\[(vg)^k=g^k.g^{-k}vg^k.g^{k-1}vg^{k-1}. \ldots. g^{-1}vg.\]
If $g$ has order $n$ then putting $k=n$ and using the fact that the successive powers of $g$ act on $V$ as the powers of $M_{\omega}$ in some order, we see from Lemma~7.1 that $(vg)^n=1$, so $vg$ has order dividing $n$. Since $vg$ maps onto the element $g$ of order $n$ under the epimorphism $G\to G/V\cong G_{\bf 0}$, it has order exactly $n$. 

Conversely, suppose that $vg$ has order $n$, so $g$ has order $k$ dividing $n$. If $k<n$ then
$1\neq (vg)^k=g^{-k}vg^k.g^{k-1}vg^{k-1}. \ldots. g^{-1}vg\in V$, so $(vg)^k$ has order $p$ and hence $k=n/p$. Thus $g$ is a primitive power of $x^p$. Since $p$ divides $n=d(q-1)$ with $q=p^e$ we see that $p$ divides $d$. However, a simple calculation, along the lines of Lemma~7.1, then shows that $\sum_{i=0}^{k-1}M_{\omega}^{pi}=0$ and hence $(vg)^k=1$, a contradiction. Thus $g$ has order $n$.\hfill$\square$

\medskip

If $d$ is even or $q$ is odd, then the order $n=d(q-1)$ of the cyclic group $G_0$ is even, so there is a unique involution $g_2\in G_0$.

\begin{cor} If $q$ is odd then the elements of order $2$ in $G$ are those of the form $vg_2$ where $v\in V$. If $q$ is even then the elements of order $2$ in $G$ are the elements of $V\setminus\{1\}$, and also, if $d=2c$ is even, those of the form $vg_2$ where $v=(v_i)\in V$ with $v_{i+c}=\sqrt\omega v_i$ for each $i=1,\ldots, c$.\hfill$\square$
\end{cor}

The proof, which we omit, is similar to that for Corollary~8.2. When $q$ and $d$ are even, the elements $v\in V$ such that $v_{i+c}=\sqrt\omega v_i$ are the fixed points of $g_2$ in $V$, forming an $F_qG$-submodule $U$ of $V$ of dimension $c=d/2$.

\medskip

If ${\cal M}\cong{\cal H}(d,\omega)$ then the standard generators $x$ and $y$ have orders $n$ and $2$, and generate $G={\rm Aut}^+{\cal M}\cong G(d,q)$. The following result shows that the converse is also true.

\begin{prop}
If elements $x$ and $y$ of orders $n$ and $2$ generate $G=G(d,q)$, then the corresponding orientably regular map is isomorphic to a Hamming map ${\cal H}(d,\omega)$ for some generator $\omega$ of $F_q^*$.
\end{prop}

\noindent{\sl Proof.} First suppose that $q$ is odd. It follows from Corollaries~8.2 and 8.3 that if $x$ and $y$ have orders $n$ and $2$ then $x^{n/2}$ and $y$ are both elements of $Vg_2$, so $y=wx^{n/2}$ for some $w\in V$ and hence $\langle x, y\rangle=\langle w, x\rangle$. Let $X=\langle x\rangle$ and let $W$ be the normal closure of $w$ in $V$, that is, the cyclic $F_pG$-submodule of $V$ generated by $w$. Then $\langle w, x\rangle=W\negthinspace :\negthinspace X$, so $W=V$ since $w$ and $z$ generate $G=V\negthinspace :\negthinspace X$. Thus $V$ is generated by $w$ as an $F_pG$-module. Now $V$ has the structure of a vector space over $F_q$, with $G$ acting linearly on $V$, so $V$ is also generated by $w$ as an $F_qG$-module. Since $G$ acts on $V$ as the cyclic group $G/V$, with the subgroup of index $d$ inducing scalar transformations, the $1$-dimensional $F_q$-subspace spanned by $w$ has at most $d$ distinct images in $V$, namely the $1$-dimensional subspaces $W_i$ spanned by $w^{x^i}$ for $i=0,\ldots, d-1$. In order for the $d$-dimensional space $V$ to be spanned by $w$ as an $F_qG$-module, these $d$ subspaces $W_i$ must be distinct and $V$ must be their direct sum; these subspaces are cyclically permuted by $X$, with $x^d$ acting on $V$ as $\omega I_d$ for some generator $\omega$ of $F_q^*$. In the orientably regular map $\cal M$ with standard generators $x$ and $y$ for $G$, each vertex is identified with a coset $gX$ of $X$ in $G$, and hence with the unique element $v \in V\cap gX$. Since $y\in wX$, one of the neighbours of the vertex $0$ is $w$, and so the neighbours of $0$ are the conjugates of $w$ in $G$. These are the non-zero elements of the subspaces $W_i$, with their cyclic ordering given by iteration of $M_{\omega}$, so ${\cal M}\cong{\cal H}(d,\omega)$. 

Now suppose that $q$ is even, so that the second part of Corollary~8.3 applies to $y$. If $y\in V$ we define $w=y$, and then the proof proceeds as before.  The other possibility is that $d$ is even and $y\in Vg_2$. But then $x^{n/2}$ and $y$ are both involutions not contained in $V$, so they are both elements of $Ug_2$, where $U$ is the subgroup of $V$ fixed by $g_2$. Thus $y=wx^{n/2}$ as before, but with $w$ contained in the proper $F_qG$-submodule $U$ of $V$. It follows that $\langle x, y\rangle=\langle w, x\rangle\leq U\negthinspace :\negthinspace X<G$, against our hypothesis that $x$ and $y$ generate $G$. This case therefore cannot arise. \hfill$\square$

\section{Galois conjugacy}

Grothendieck's theory of {\it dessins d'enfants}~\cite{Gro, JS2} shows that a map $\cal M$ on a compact oriented surface corresponds naturally to a {\sl Bely\u{\i} pair\/} $(X,\beta)$, where $X$ is a nonsingular projective algebraic curve over $\bf C$, and $\beta$ is a rational function from $X$ to the complex projective line (or Riemann sphere) ${\bf P}^1({\bf C})={\bf C}\cup\{\infty\}$, unramified outside $\{0, 1, \infty\}$. One can regard $X$ as a Riemann surface underlying $\cal M$, with the inverse image under $\beta$ of the unit interval providing the embedded graph. By Bely\u{\i}'s Theorem~\cite{Bel}, $X$ and $\beta$ are defined (by polynomials and rational functions) over the field $\overline{\bf Q}$ of algebraic numbers, and conversely every curve defined over $\overline{\bf Q}$ is obtained in this way from a map. Grothendieck observed that the action of the {\sl absolute Galois group\/} ${\bf \Gamma}={\rm Gal}\,\overline{\bf Q}$ on the coefficients of these functions induces a faithful action of $\bf \Gamma$ on the associated maps. It is therefore of interest to determine the orbits of $\bf \Gamma$ on maps.

The following result extends examples given by Streit, Wolfart and the author in~\cite{JSW}; in particular it generalises the corresponding result~\cite[Theorem 4]{JSW} for complete graphs $H(1,q)=K_q$:

\begin{thm} For any $d\geq 1$ and prime power $q=p^e$, the $\phi(q-1)/e$ Hamming maps ${\cal H}(d,\omega)$ form an orbit under $\bf \Gamma$.
\end{thm}

\noindent{\sl Proof.} As shown by Streit and the author~\cite{JSt}, various properties of a map $\cal M$ are invariant under the action of $\bf \Gamma$: these include orientable regularity, the (orientation-preserving) automorphism group ${\rm Aut}^+{\cal M}$, and the type $\{m,n\}$ of $\cal M$. Proposition~8.4 shows that for any given $d$ and $q$ the Hamming maps ${\cal H}(d,\omega)$ are characterised among all orientably regular maps by their common automorphism group $G(d,q)$ and their valency $n=d(q-1)$, so this set of maps is invariant under $\bf \Gamma$. Lemma~3.2 shows that these maps form a single orbit under Wilson's operations $H_j$. Now Streit, Wolfart and the author have shown in~\cite[Theorem~2]{JSW} that any set of maps satisfying these two conditions forms an orbit of $\Gamma$, so in particular this applies to the Hamming maps ${\cal H}(d,\omega)$. \hfill$\square$

\medskip

Theorem~2 of~\cite{JSW} also allows one to determine the minimal field of definition of the Bely\u{\i} pairs corresponding to these maps. For any integer $m\geq 1$ let $\zeta_m=\exp(2\pi i/m)\in{\bf C}$. The multiplicative group ${\bf Z}^*_m$ of units mod~$(m)$ can be identified with the Galois group of the cyclotomic field ${\bf Q}(\zeta_m)$ by identifying each $j\in {\bf Z}^*_m$ with the automorphism defined by $\zeta_m\mapsto\zeta_m^j$. It follows from~\cite[Theorem 2]{JSW} that the Bely\u{\i} pairs corresponding to the maps ${\cal H}(d,\omega)$ are all defined over a particular subfield of ${\bf Q}(\zeta_n)$, namely the fixed field $K$ of the subgroup $H$ of ${\rm Gal}\,{\bf Q}(\zeta_n)$ consisting of all $j\in {\bf Z}^*_n$ such that $H_j({\cal H}(d,\omega))\cong{\cal H}(d,\omega)$ for some (and hence every) $\omega$. By Theorem~2.2 and Lemma~3.2, $H$ is the inverse image of the subgroup $\langle p\rangle\leq{\bf Z}_{q-1}^*$ under the natural epimorphism ${\bf Z}_n^*\to{\bf Z}_{q-1}^*$. This is a subgroup of index $\phi(q-1)/e$ in ${\bf Z}_n^*$, containing the kernel of this epimorphism, so $K$ is an extemsion of $\bf Q$ of degree $\phi(q-1)/e$, contained in ${\bf Q}(\zeta_{q-1})$, namely the subfield of ${\bf Q}(\zeta_{q-1})$ fixed by its automorphism $\zeta_{q-1}\mapsto \zeta_{q-1}^p$. Equivalently, $K$ is the splitting field of $p$ in ${\bf Q}(\zeta_{q-1})$, the maximal subfield of ${\bf Q}(\zeta_{q-1})$ in which $p$ decomposes into $\phi(q-1)/e$ different primes ideals of degree $1$ over $p$.

\medskip

\noindent{\bf Example 1.} The Bely\u{\i} pair corresponding to ${\cal H}(d,\omega)$ is defined over $\bf Q$ if and only if $\phi(q-1)/e=1$, that is, $q=2, 3$ or $4$. As a simple example, the unique orientably regular embedding of $H(2,2)\cong Q_2\cong K_{2,2}$ corresponds to the Fermat curve $X$ of genus $0$ given (as a projective curve) by $x^2+y^2=z^2$, with Bely\u{\i} function $\beta:[x,y,z]\mapsto 4x^2(z^2-x^2)/z^4$. The vertices are the points $[0,\pm 1, 1]$ and $[\pm 1,0,1]$ where $\beta=0$, and the edges are the points where $0\leq\beta\leq 1$, with their centres at the points $[\pm 1/\sqrt 2, \pm 1/\sqrt 2, 1]$ where $\beta=1$.

\medskip

\noindent{\bf Example 2.} If $q=25$ then for each $d\geq 1$ there are four maps ${\cal H}(d,\omega)$, described in \S 4(b). The corresponding Bely\u{\i} pairs are defined over the splitting field $K$ of the prime $p=5$ in the cyclotomic field ${\bf Q}(\zeta_{24})$. This is an extension of $\bf Q$ of degree $4$, with the four maps forming an orbit under ${\rm Gal}\,K$ (a Klein four-group generated by $\zeta_{24}\mapsto\zeta_{24}^{-1}$ and $\zeta_{24}\mapsto\zeta_{24}^7$), and hence also forming an orbit under $\bf \Gamma$, which acts as ${\rm Gal}\,K$ via the epimorphism ${\bf \Gamma}\to{\rm Gal}\,K$ induced by the inclusion of $K$ in $\overline{\bf Q}$.

\medskip

In general, finding explicit equations for a Bely\u{\i} pair, as in Example~1, is much harder than finding its Galois orbit or minimal field of definition. It would be interesting to achieve this for some more of the Hamming maps.

\section{Distance $k$ Hamming graphs}

The distance $k$ Hamming graphs $H_k=H(d,q)_k$, for $k=1, 2, \ldots, d$, are generalisations of the Hamming graphs $H=H(d,q)$. They have the same vertex set $V=Q^d$ as $H$, with vertices $v=(v_i)$ and $w=(w_i)$ adjacent in $H_k$ if and only if they are at distance $k$ in $H$, that is, $v_i\neq w_i$ for exactly $k$ values of $i$. Thus $H_1=H$. The vertices of $H_k$ all have valency ${d\choose k}(q-1)^k$. These graphs are arc-transitive, since $A={\rm Aut}\,H=S_q\wr S_d$ is a group of automorphisms of $H_k$, acting transitively on its arcs.

For any subset $K$ of $D:=\{1, 2, \ldots, d\}$ let $H_K=H(d,q)_K$ be the merged Hamming graph
$\cup_{k\in K}H(d,q)_k$: this has vertex set $V$, with $v$ and $w$ adjacent if and only if they are at Hamming distance $k$ for some $k\in K$. We have $A\leq{\rm Aut}\,H(d,q)_K$ for each $K$, and when $q>3$ we can determine whether these two groups are equal. (We will explain later why the case $q=3$  is excluded here.) Let $D_0$ and $D_1$ denote the sets of even and odd elements of $D$.

\begin{prop}
Let $q\geq 4$ and let $K$ be a non-empty subset of $D=\{1, 2, \ldots, d\}$. Then ${\rm Aut}\,H(d,q)_K=A$ in all cases except the following:
\begin{description}
\item[(a)] $K=D$;
\item[(b)] $q=4$, $d\geq 3$ and $K=D_0$ or $D_1$.
\end{description}
\end{prop}

\noindent{\sl Proof.} As the automorphism group of a graph, $A$ is a $2$-closed permutation group, that is, it is the full automorphism group of the set of binary relations on $V$ which it preserves. Any permutation group $A^*$ on $V$, which properly contains $A$, must therefore have rank less than the rank $d+1$ of $A$, so the association scheme ${\cal A}^*$ corresponding to $A^*$ must be a proper subscheme of the Hamming association scheme $\cal A$ corresponding to $A$, which has the adjacency matrices $A_k$ of the graphs $H_k\;(k=0,\ldots, d)$ as its basis elements. Muzychuk~\cite{Muz} has shown that for $q\geq 4$ the only such subschemes are the rank~$2$ association scheme with basis $\{A_0=I, A_1+A_2+\cdots+A_d\}$, and the rank~$3$ association scheme for $q=4$ and $d\geq 3$ with basis $\{A_0=I, A_1+A_3+\cdots, A_2+A_4+\cdots\}$. If ${\rm Aut}\,H(d,q)_K\neq A$ then by applying this to $A^*={\rm Aut}\,H(d,q)_K$ we see that $K$ is as described in (a) or (b). \hfill$\square$

\medskip

For $q\geq 5$, the fact that ${\rm Aut}\,H(d,q)_K=A$ can also be deduced from the maximality of the group $A=S_q\wr S_d$ in the alternating or symmetric group on $V$: there are proofs of this by Neumann in~\cite{JSo} and by Liebeck, Praeger and Saxl in~\cite{LPS}, both depending on the classification of finite simple groups, and there is a more elementary proof, valid for $q$ sufficiently large as a function of $d$, by Soomro and the author in~\cite{JSo}. Muzychuk's classification of subschemes of the Hamming schemes for $q\geq 4$ is purely combinatorial, and is also independent of the classification of finite simple groups.

In the exceptional cases in Proposition~10.1, if $K=D$ then $H_K$ is a complete graph and so its automorphism group is the symmetric group $S_{q^d}$. If $q=4$ then the two complementary graphs $H_K$ corresponding to $K=D_0$ and $D_1$ have the same rank $3$ automorphism group. This can be explained as follows.

Let $\omega\in F_4\setminus F_2$. Given any $v=(v_i)\in V=F_4^d$, one can write each $v_i$ uniquely as $x_i+\omega y_i$ with $x_i, y_i\in F_2$, thus identifying $V$ with $F_2^{2d}$. Now $x_i^2+x_iy_i+y_i^2=0$ or $1$ as $v_i=0$ or $v_i\neq 0$, so the subgroup $A_0=S_3\wr S_d$ of $A$ fixing $0$, acting on $V$ as $GL_2(2)\wr S_d$, preserves the quadratic form
\[Q(v)=\sum_{i=1}^d(x_i^2+x_iy_i+y_i^2),\]
and is therefore a subgroup of the orthogonal group corresponding to $Q$. This is the group $GO^{\varepsilon}_{2d}(2)$ in ATLAS notation~\cite{CCNPW}, where the sign $\varepsilon$ is $+$ or $-$ as $Q$ has Witt index $d$ or $d-1$, that is, as $d$ is even or odd. It follows that the semidirect product $A^*=V\negthinspace:\negthinspace GO^{\varepsilon}_{2d}(2)$, with $GO^{\varepsilon}_{2d}(2)$ acting naturally on $V$, is a subgroup of the affine group $AGL_{2d}(2)$ containing $A=V\negthinspace:\negthinspace A_0$. Now $GO^{\varepsilon}_{2d}(2)$ has two orbits $\Gamma_0$ and $\Gamma_1$ on $V\setminus\{0\}$, consisting of the isotropic and non-isotropic vectors, those $v\neq 0$ with $Q(v)=0$ or $1$, so $A^*$ acts on $V$ as a rank $3$ group. For $K=D_0$ or $D_1$, the complementary graphs $H_K$ on $V$ are defined by $v$ being adjacent to $w$ if and only if $v-w$ has even or odd Hamming weight (distance from $0$), that is, $v-w\in\Gamma_0$ or $\Gamma_1$ respectively. These relations are invariant under $V$ and $GO^{\varepsilon}_{2d}(2)$, so in each case ${\rm Aut}\,H_K$ contains $A^*$ and therefore has rank $3$.

\begin{thm}
Let $q\geq 4$ and let $K$ be a non-empty subset of $D=\{1, 2, \ldots, d\}$. Then the graph $H_K=H(d,q)_K$ has an orientably regular embedding if and only if $q$ is a prime power and $K=\{1\}$ or $D$.
\end{thm}

\noindent{\sl Proof.} Let $q$ be a prime power. If $K=\{1\}$ then $H_K=H$ and so the existence of an orientably regular embedding follows from Theorem~2.1. If $K=D$ then $H_K$ is a complete graph, and the number $q^d$ of vertices is a prime power, so the embeddings constructed by Biggs~\cite{Big} give the result.

For the converse, suppose that $\cal M$ is an orientably regular embedding of $H_K$, and let $G={\rm Aut}^+{\cal M}$. Then ${\rm Aut}\,H_K$ has a cyclic subgroup $G_0$ fixing $0$ and acting regularly on its neighbours, that is, on the vectors of weight $k\in K$. First suppose that $K$ is not one of the exceptional subsets in Proposition~10.1, so that ${\rm Aut}\,H_K={\rm Aut}\,H=A$. Then $G_0\leq A_0$, so $G_0$ preserves the weights of all vectors and hence $K$ is a singleton $\{k\}$ for some $k\in D$. 
Since $G_0$ acts transitively on the vectors of weight $k$ it must map onto a subgroup of $A/B\cong S_d$ acting transitively on the $k$-element subsets of $D$; since $G_0$ is abelian, this subgroup acts regularly on $D$ and hence $k=1$, $d-1$ or $d$. 

If $k=1$ then $H_K=H$, so Theorem~2.1 implies that $q$ is a prime power.

Next suppose that $k=d-1$. Then $G_0$, a cyclic group of order $d(q-1)^{d-1}$, must map onto a cyclic subgroup of $S_d$ acting regularly on $D$. The kernel $G_0\cap B$, a cyclic subgroup of $B_0=S_{q-1}^d$ of order $(q-1)^{d-1}$, must map onto a transitive and hence regular subgroup of each factor $S_{q-1}$. The unique subgroup of index $q-1$ in $G_0\cap B$ must therefore be trivial, so $d=2$ and hence $k=1$, a case we have already dealt with.

Now let $k=d$, so $G_0$ is a cyclic subgroup of $A_0=S_{q-1}\wr S_d$, acting regularly on the set of vectors of weight $d$. Since $q\geq 4$, $A_0$ acts primitively on this set. Since $A_0$ has a cyclic regular subgroup, and the degree $(q-1)^d$ is not prime, a theorem of Burnside~\cite[\S 252]{Bur} and Schur~\cite{Sch} implies that $A_0$ must act doubly transitively, whereas in fact it has rank $d+1>2$. Thus there are no orientably regular embeddings in this case.

Now suppose that $K$ is one of the exceptional sets in Proposition~10.1, so that we have ${\rm Aut}\,H_K>A$. If $K=D$ then $H_K$ is a complete graph on $q^d$ vertices. As shown by Biggs~\cite{Big}, if this has an orientably regular embedding then $q^d$ is a prime power, and hence so is $q$. We may therefore assume that $q=4$, $d\geq 3$ and $D=D_i$ for $i=0$ or $1$. By the comments following Proposition~10.1, ${\rm Aut}\,H_K$ contains the group $A^*=V\negthinspace:\negthinspace GO^{\varepsilon}_{2d}(2)$. Now $A^*_0=GO^{\varepsilon}_{2d}(2)$ acts primitively on each of its orbits $\Gamma_0$ and $\Gamma_1$ in $V\setminus\{0\}$, and $G_0$ induces a cyclic permutation on the set $\Gamma_i$ of neighbours of $0$ in $H_K$. It follows from the theorem of Burnside and Schur mentioned above that $({\rm Aut}\,H_K)_0$, which contains both of these groups and does not act as a subgroup of $AGL_1(p)$ for any prime $p$, must act as a doubly transitive group on $\Gamma_i$. However, this is impossible, since in either graph $H_K$ it is easy to find pairs of neighbours of $0$ which are and are not adjacent.\hfill$\square$

\medskip

If $q=3$ then the preceding arguments in the case ${\rm Aut}\,H_K=A$ yield one further example with an orientably regular embedding, namely $H(2,3)_2$. Once again, $K$ must be a singleton $\{k\}$, with $k=1$, $d-1$ or $d$, and the first two cases give the same conclusions as before. In the case $k=d$, however, $A_0=S_{q-1}\wr S_d$ is now imprimitive on vectors of weight $d$, so the preceding argument does not apply. Instead, note that a cyclic regular subgroup $G_0$ of order $2^d$ in $A_0$ would map onto a cyclic group of order $2^{d-1}$ in $S_d$, impossible for $d>2$ since a generator would have to contain a cycle of length $2^{d-1}>d$. Hence $d=2$, so $H_K=H(2,3)_2$, which is isomorphic to $H(2,3)$ and therefore has a unique orientably regular embedding. To deal with cases where ${\rm Aut}\,H_K>A$ one would need to know the subschemes of the Hamming scheme for $q=3$. One of these arises from a rank $4$ group $A^*=V\negthinspace :\negthinspace GO_d(3)\geq A$, where the orthogonal group preserves the quadratic form $\sum_{i=1}^dv_i^2$, which takes the value $0, 1$ or $2$ as $v$ has weight congruent to $0, 1$ or $2$ mod~$(3)$. However, at present there is no complete analogue for $q=3$ of Muzychuk's classification of subschemes of the Hamming scheme for $q\geq 4$. It is hoped to consider this problem in a future paper.

\end{document}